# The strange story of an almost unknown prime number counter: The Rafael Barrett formula

Eduardo Mizraji

Group of Cognitive Systems Modeling, Biophysics and Systems Biology Section,
Faculty of Sciences, University of the Republic, Uruguay
ORCID ID: https://orcid.org/0000-0001-6938-8427
mizraj@fcien.edu.uy; emizraji@gmail.com

ABSTRACT

In this brief article, we present the formula created by Rafael Barrett in 1903 in a note to Henri Poincaré, which remained unknown for decades. Discovered in the 1930s by a Uruguayan mathematician, this formula was published and analyzed in a journal published in Montevideo in 1935. In this study, we present Barrett's formula and analyze a challenge it could pose.

In this short article, we want to share an interesting and educational story. It is about how an eminent writer named Rafael Barrett, who had a great influence on the literature and culture of Spanish-speaking countries, was also a secret lover of mathematics. Many years after his death, a Uruguayan mathematician found a formula created by Barrett that counted prime numbers $p < n$, where n is any natural number. Barrett had intended to send it to Henry Poincaré, and the brief letter to the French mathematician containing the formula is what the Uruguayan mathematician found. As we shall see, Barrett's formula is not only ingenious, it is also beautiful. We wish to emphasize here the overall significance of this situation, in which a prominent writer possesses a secret world where he explores challenging territories of mathematics. This is evidence that there are fortunate minds in which several worlds can coexist.

Rafael Barrett was a writer and essayist, a "*homme de lettres",* gifted with virtuoso writing (Figure 1). He was born in Spain in 1876, the son of an Englishman and a Spanish woman. After starting his university studies, he moved to South America, where, after passing through Argentina, he settled in Asunción, Paraguay. Barrett was moved by the tragedy of this country, which, after being destined to be one of the economic and cultural powers of South America, had been devastated a few decades earlier in a war against Brazil, Argentina, and Uruguay, which plunged it into extreme poverty (the Paraguayan War, or War of the Triple Alliance, took place between 1864 and 1870).

After his stay in Asunción, Barrett moved to Buenos Aires and Montevideo, where these cities of Argentina and Uruguay had an intense cultural life in which Barrett's talent was immediately detected by the important literary figures of those years. A good outline of his biography in English can be found in Wikipedia (https://en.wikipedia.org/wiki/Rafael_Barrett). A description of the life of Barrett in South America can be found in [1] and a detailed analysis of his work and the cultural context of his generation was published in Spanish by Czech researcher Maksymilian Drozdowicz of the University of Ostrava [2].

.

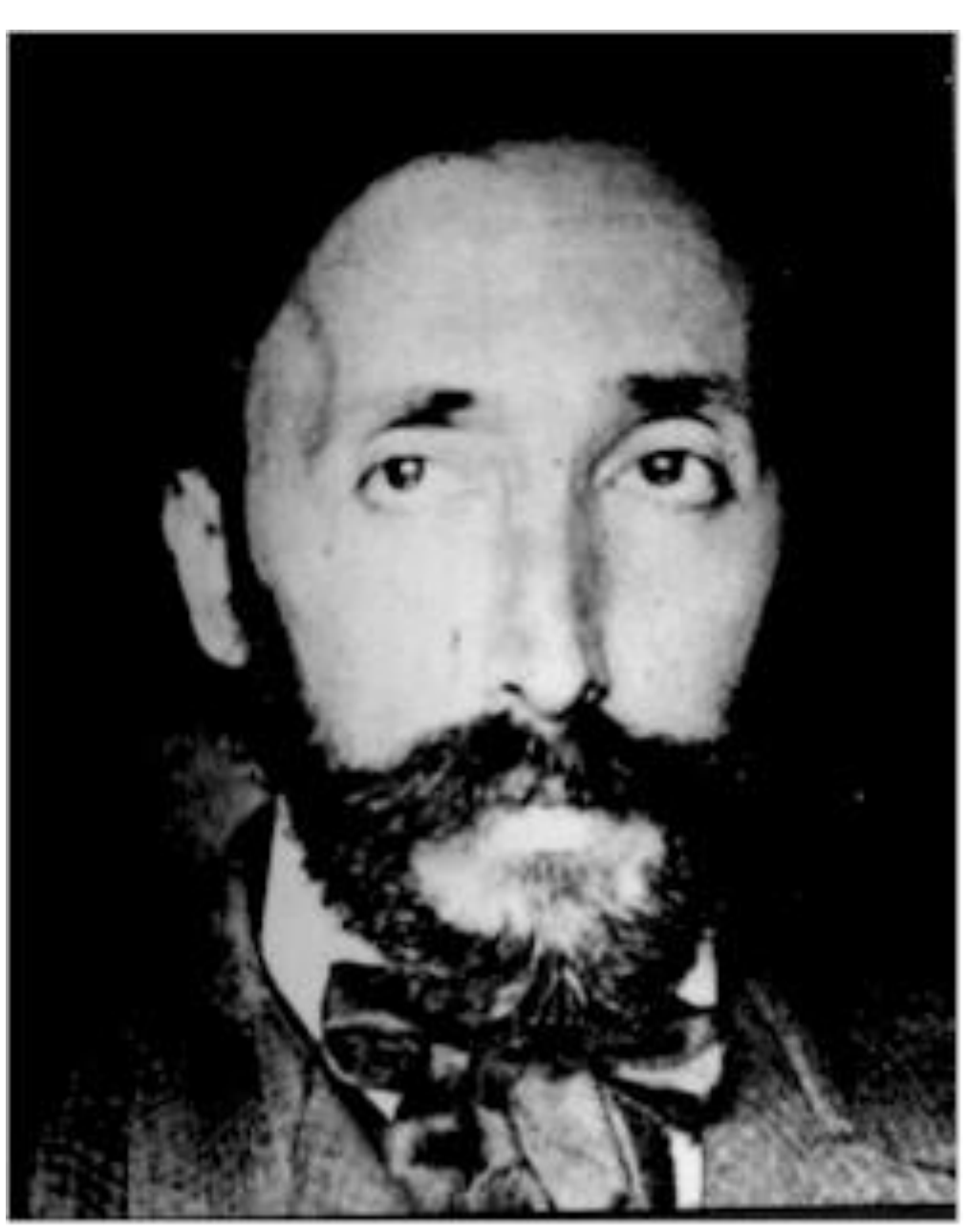

Fig. 1 Rafael Barrett (public domain image from the respective article in Wikipedia)

In 1910, Barrett died in France at the age of 34 from a lung condition. Posthumously, a beautiful book of essays by Barrett, "Mirando Vivir" [3], was published in Montevideo, which contains texts that cover literary, political and scientific topics and testifies to the breadth of Barrett's interests and thought. *"Paraguayan Sorrow: Writings of Rafael Barrett*, a collection of essays originally published in Montevideo in 1911, is the only one of his writings that we know has been translated into English [4].

Here the focus is another side of Barrett, which does not seem to have been documented until 1935. In this year, Eduardo García de Zúñiga published an article in Montevideo in which he analyzes a formula that Barrett wrote in 1903 to Henri Poincaré containing a prime number counter created by him [5]. Eduardo García de Zúñiga was the "father" of the Uruguayan mathematical community, whose outstanding research contributions are evident in the work of Massera and Schäffer on differential equations [5] and whose research later expanded to the present day to stochastic processes and algebra, among other topics. Figure 2 shows Professor García de Zúñiga.

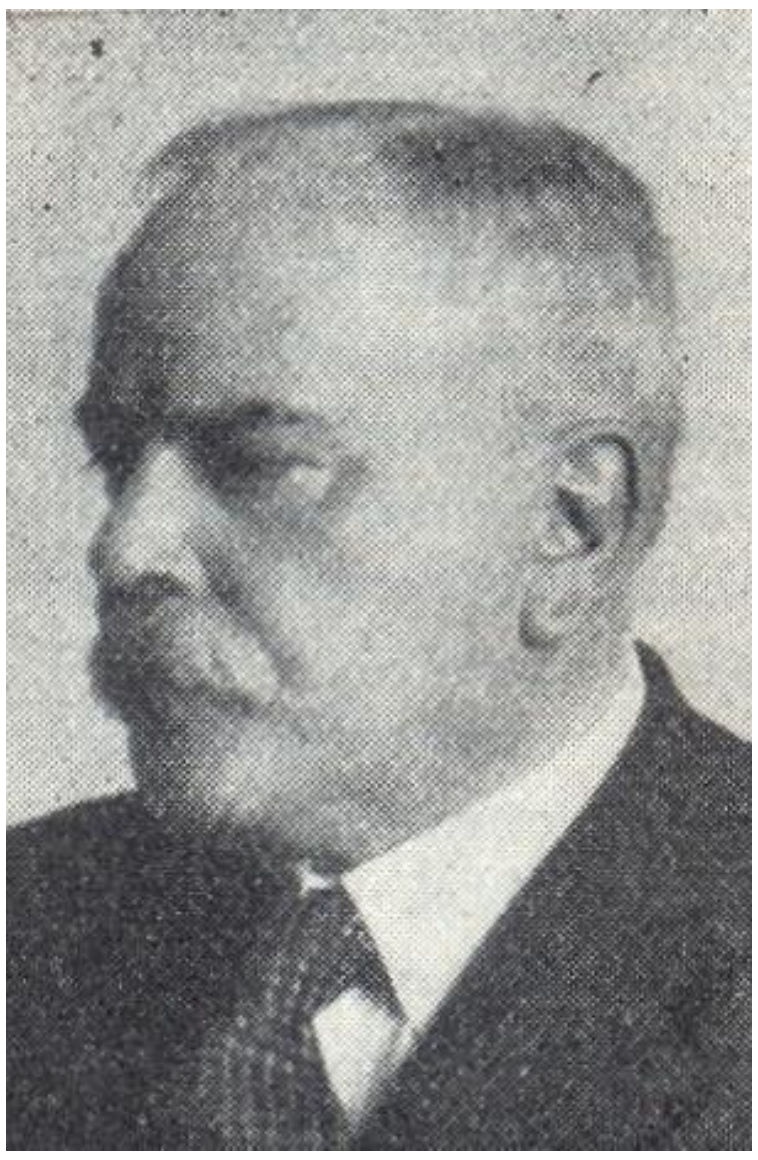

Fig. 2. Eduardo García de Zúñiga (public domain image from the respective article in Wikipedia)

Figure 3 shows what Rafael Barrett prepared for Poincaré taken from the article by García de Zúñiga (who in 1935 was not sure if this proof had actually been sent).

*A M. H. Poincaré, de l'Académie des Sciences.*

*Monsieur,*

*J'ai l'honneur de vous communiquer une formule relative au nombre* $\lambda$ *des nombres premiers inférieurs à une limite donnée* L.

*La fraction*

$$\frac{1.\ 2.\ 3.\ldots..\ (n-1)}{n} \quad (n \text{ entier} > 4)$$

*est égale à un nombre entier pour n composé, et égale à un nombre entier moins 1/n pour n premier.*

*Donc,*

$$\sin \frac{1.\ 2.\ 3.\ldots..\ (n-1)\ \pi}{n} \quad (n \text{ entier} > 4)$$

*sera égal à 0 ou à 1, suivant que n sera composé ou premier.*

*Alors, si V est le nombre entier immédiatement inférieur à la limite donnée* L, *le nombre* $\lambda$ *des nombres premiers inférieurs à* L *sera donné par la série*

$$3 + \sum_{n=5}^{n=1} \frac{\sin \dfrac{(n-1)!\ \pi}{n}}{\sin \dfrac{\pi}{n}}$$

RAFAEL BARRETT.

*6 Oct., 1903.*

Fig. 3. Text of the note written by Barrett to Poincaré (with a small correction by García de Zúñiga).

The problem Barrett addresses is a classic theme of prime number theory. Gauss and Lagrange conjectured that for a natural number n, the asymptotic distribution of prime numbers $p < n$ is given by the formula:

$$\Pi(n) \sim \frac{n}{\ln(n)} \quad , \quad n \to \infty \ ; \qquad (1)$$

This conjecture was proven by Hadamard and La Vallée Poissin in 1896 [6].

The following is Barrett's formula (we will note Barr(n)), which expresses the number of prime numbers $p < n$:

$$\mathrm{Barr}(n) = 3 + \sum_{k=5}^{n-1} \frac{\sin\left[\dfrac{(k-1)!}{k}\pi\right]}{\sin\left(\dfrac{\pi}{k}\right)} \ . \qquad (2)$$

Let us now show, with a miniature example, how this formula behaves:

Barr(6) = Barr(7) = 4; Barr(8) = Barr(9) = Barr(10) = Barr(11) = 5; Barr(12) = Barr(13) = 6; Barr(14) = Barr(15) = Barr(16) =Barr(17) = 7; etc. (At that time, Barrett included the number 1 among the primes.)

García de Zúñiga showed that this formula arises from Wilson's theorem on prime numbers [5, 7]:

$$\text{p prime} \Leftrightarrow (p-1)!+1 \equiv 0 \pmod{p} \ . \qquad (3)$$

In [5], there is a detailed proof of Wilson's theorem that shows how Barrett arrived at his formula. We include a proof of Barrett's formula in the Appendix.

It is clear that Barrett's formula cannot be associated with the computational efficiency of a prime number counter. For the evolution and current references on prime counters and computational efficiency, see [8].

It is interesting to consider that an unbounded prime counter implicitly implies the limit (1). Accessing that limit from a prime counter might be impossible. But if this impossibility is not proven, we could conjecture that perhaps that limit is reachable if someone finds an ingenious heuristic. In this regard, an eloquent illustration of a somewhat surprising heuristic is the procedure shown by Courant and Robbins in the last section of their celebrated book “What is Mathematics?” [9, page 482]. The title of this section is “The prime number theorem obtained by statistical methods”. The authors state that the procedure was suggested to them by Gustav Hertz (an experimental physicist who obtained the Nobel Prize in 1925). They begin with two approximations: Stirling’s formula expressed as $\ln(n!) \sim n\ln(n)$ and Legendre’s theorem (shown in Appendix, eq. A2) for n!, without the floor function and converted into a geometric series of primes. From there, in a few steps and further approximations, they elegantly arrive at the limit (1).

In the Appendix that follows, we present a detailed and elementary proof of Barrett's formula, omitting no steps or condensing reasoning. We also take the opportunity to illustrate and utilize Legendre's important theorem, which, incidentally, provides the reader with a tool for understanding Courant and Robbins' ingenious approach to the prime number theorem.

This concludes the report we intended to produce, which demonstrates that sometimes the talent dissolves the confrontation between diverse cultural modes. As a final point, and assuming that any prime number counter implies the prime number theorem, however difficult or improbable what follows may seem, perhaps we should ask: *Can there be any heuristic that can obtain the limit (1) of Barrett's formula?*

Acknowledgments. I thank the reviewers and the editor for their constructive suggestions.
Conflict of interest: The author has no interest to declare.
Funding: No funding was received

APPENDIX: Proof of Barrett's formula

Here is the proof of Barrett's formula. We take the terms of the summation,

$$\frac{\sin\left[\frac{(k-1)!}{k}\pi\right]}{\sin\left(\frac{\pi}{k}\right)} \qquad (4)$$

and we analyze what happens if k is composite or if k is prime.

Previously, let us recall that $\sin(n\pi)=0$ and $\cos(n\pi)=(-1)^n$ for all integers n. We now analyze two possible situations.

1) k is composite

In this case, the prime factors of k are in $(k-1)!$. Here follows the argument that proves this assertion. Let $k=\Pi p_i^{\alpha_i}$. Legendre's formula allows us to determine the exponents of the prime factors of a factorial

$$v_{p_i}((k-1)!)=\sum_{j\geq 1}\left\lfloor \frac{k-1}{p_i^{\,j}} \right\rfloor \;, \qquad (5)$$

where each term is the integer part ("floor" function) of the displayed quotients. From this theorem we get that when $k > 4$ then k is a factor of $(k-1)!$. Hence k divides $(k-1)!$ and the numerator of (4) is $\sin(n\pi)=0$.

2) k is prime

In this case, we analyze the quotient (A3). We write it as follows:

$$\frac{(k-1)!+1-1}{k}=\frac{[(k-1)!+1]}{k}-\frac{1}{k} \;. \qquad (6)$$

But by Wilson's theorem, since k is prime, $(k-1)!+1\equiv 0(\text{mod}\, k)$ is an odd integer. Then, standard properties of the sine show that $\sin[((k-1)!)\pi)/k]=\sin(\pi/k)$. Therefore, if k is prime the quotient (4) gives 1.

This explains the Barrett equation as a prime counter $p < n$. □